\documentclass[12pt,twoside]{article}
\usepackage{latexsym}
\usepackage{amssymb}
\usepackage{amsfonts}
\makeatletter
\newtheorem{theorem}{Theorem}[section]
\newtheorem{lemma}[theorem]{Lemma}

\newtheorem{defin}[theorem]{Definition}
\newtheorem{result}[theorem]{Result}

\def\ps@headings{
 \def\@oddhead{\footnotesize\rm\hfill\runningheadodd\hfill\thepage}
 \def\@evenhead{\footnotesize\rm\thepage\hfill\runningheadeven\hfill}
 \def\@oddfoot{} \def\@evenfoot{\@oddfoot}
}

\newcommand{\Prf}{\noindent{\bf Proof}.\quad }

\newcommand{\qed}{\hfill$\Box$}

\makeatother

\def\runningheadeven{TRANSITIVE AND CO--TRANSITIVE CAPS}
\def\runningheadodd{A. Cossidente and O.H. King}

\title{Transitive and Co--Transitive caps}

\author{{\em A. Cossidente} , {\em O.H. King}}

\begin{document}
\maketitle
\pagestyle{headings}

\section{Introduction}

Let $PG(r,q)$ be the projective space of dimension $r$ over
$GF(q)$. A $k$--{\it cap} $\bar{K}$ in $PG(r,q)$ is a set of $k$ points,
no three of which are collinear \cite{HT}, and a $k$--cap is said to be
{\it complete} if it is maximal with respect to set--theoretic inclusion.
The maximum value of $k$ for which there is known to exist a $k$--cap in
$PG(r,q)$ is denoted by $m_2(r,q)$. Some known bounds for $m_2(r,q)$ are given  below. 

Suppose that $\bar{K}$ is a cap  in $PG(r,q)$ with automorphism group $\bar{G}_0 \leq
P\Gamma L(r+1,q)$. Then $\bar{K}$ is said to be  {\em transitive} if $\bar{G}_0$ acts transitively on $\bar{K}$, and {\em co-transitive} if $\bar{G}_0$ acts transitively on $PG(r,q) - \bar{K}$.

Our main result is the following theorem.

        \begin{theorem}
Suppose $\bar K$ is a transitive, co--transitive cap in $PG(r,q)$.
Then one of the following occurs:
\begin{enumerate}

\item $\bar K$ is an elliptic quadric in $PG(3,q)$ and $q$ is a
square when $q$ is odd;
\item $\bar K$ is the Suzuki--Tits ovoid in $PG(3,q)$ and $q=2^h$,
with $h$ odd and $\geq 3$;
\item $\bar K$ is a hyperoval in $PG(2,4)$;
\item $\bar K$ is an $11$--cap in $PG(4,3)$ and
$\bar{G_0}\simeq M_{11}$;
\item $\bar K$ is the complement of a hyperplane in $PG(r,2)$;
\item $\bar K$ is a union of Singer orbits in $PG(r,q)$ and
$G_0\leq\Gamma L(1,p^d)\leq GL(d,p)$.
\end{enumerate}
In each of $1$--$5$  $\bar K$ is indeed a transitive
co--transitive cap.

        \end{theorem}

Our conclusion is that transitive, co--transitive caps are rare
with the possible exception of unions of Singer cyclic orbits.

The origin of this problem are papers by Hill \cite{RH2},
\cite{RH1}, in which he studies such caps whose automorphism group
acts $2$--transitively on the cap. [As he notes \cite[Theorem 1]{RH2}, it
is trivial to show that if $\bar{K}$ is a subset of $PG(r,q)$ lying in no
proper subspace and admitting a $3$--transitive group then $\bar{K}$ must
be a cap.] Hill gives a short list of possibilities (omitting Suzuki--Tits
ovoids) but excludes caps in $PG(r,q)$ for $q>2$ and $r\geq 13$. We find
no new caps but show that any other transitive, co--transitive cap is a
union of Singer cyclic orbits.

The known upper bounds on cap sizes are summarised in the following Result.

\begin{result}\cite[Theorem 27.3.1]{HT} \label{bound}\\
$m_2(2,q)=q+1$ (for $q$ odd);\\
$m_2(2,q)=q+2$ (for $q$ even);\\
$m_2(3,q)=q^2+1$ for $q>2$; \\
$m_2(r,2)=2^r$; and\\
$m_2(r,q)\leq q^{r-1}$ for $q>2$ and $r\geq 4$.
\end{result}

We begin by showing that as a consequence of Result \ref{bound}, a cap must be smaller than its complement (with one exception). It then follows that in considering subgroups of $P\Gamma L(r+1,q)$ having two orbits, we need only consider the smaller orbit when looking for transitive, co-transitive caps.

\begin{lemma}
Suppose that $\bar{K}$ is a cap in $PG(r,q)$. Then either $| \bar{K} | < (q^{r+1}-1)/2(q-1)$, or $q=2$ and $\bar{K}$ is the complement of a hyperplane.
\end{lemma}

\Prf It is easy to deduce from Result \ref{bound}, that the result holds when $q \neq 2$. Thus suppose now that $q=2$ and that $| \bar{K} | \geq (2^{r+1}-1)/2$. The only possiblity is that $| \bar{K} | = 2^2$. Let $x \in \bar{K}$. For each $y \in \bar{K} - \{x\}$ there is a line through $x$ and $y$ and the $2^r-1$ such lines must be distinct since $\bar{K}$ is a cap. However $x$ lies on exactly $2^r-1$ lines in $PG(r,2)$ and so every line in $PG(r,2)$ through $x$ meets $\bar{K}$ in two points and $PG(r,q) - \bar{K}$ in one point. Therefore any line meeting $PG(r,q) - \bar{K}$ in at least two points is contained in $PG(r,q) - \bar{K}$. This shows that $PG(r,q) - \bar{K}$ is a subspace of $PG(r,2)$ and its size shows that it is a hyperplane.

\qed
\\
Using Result \ref{liebeck} we shall know orbit lengths when looking at candidates for transitive, co-transitive caps. Lemma \ref{chordno} below helps in eliminating a number of possibilities.

\begin{defin}
Suppose that $\bar{K}$ is a cap in $PG(r,q)$. For any $x \in PG(r,q)$, the {\em chord-number} of $x$ is the number of chords (2-secants) of $\bar{K}$ passing through $x$.
\end{defin}

\begin{lemma}\label{chordno}
Suppose that $\bar{K}$ is a tranistive, co-transitive cap in $PG(r,q)$ and suppose that $x \in PG(r,q) - \bar{K}$. Let $k=|\bar{K}|$ and $m=|PG(r,q)-\bar{K}|$. Then the chord-number, $c$, of $x$ is given by
\begin{center}
$c={k(k-1)(q-1) \over 2m}$.
\end{center}
In particular the expression for $c$ always yields an integer.
\end{lemma}

\Prf We count combinations of chords and points of $PG(r,q) - \bar{K}$ in two ways. Firstly there are $k(k-1)/2$ chords of $\bar{K}$ and each has $q-1$ points not in $\bar{K}$. There is a subgroup $\bar{G_0}$ of $\Gamma L(r+1,q)$ acting transitively on $PG(r,q) - \bar{K}$, so each of these $m$ points has the same chord-number $c$ and a second count gives $mc$ chord-point combinations. Thus $mc=k(k-1)(q-1)/2$ leading to the required expression for $c$.

\qed

The main tool in our investigation is the substantial result by M.W. Liebeck \cite{L}, where the affine permutation groups of rank three are classified.

        \begin{result}\cite{L} \label{liebeck}
Let $G$ be a finite primitive affine permutation group of rank
three and of degree $n=p^d$, with socle $V$, where $V\simeq
(Z^d_p)$ for some prime $p$, and let $G_0$ be the stabilizer of
the zero vector in $V$. Then $G_0$ belongs to one of the following
families:
\par\noindent
{\bf(A)} 11 Infinite classes;
\par\noindent
{\bf (B)} Extraspecial classes with $G_0\leq N_{\Gamma L(d,p)}(R)$,
where $R$ is a $2$--group or $3$--group irreducible on $V$;
\par\noindent
{\bf(C)} Exceptional classes. Here the socle $L$ of $G_0/Z(G_0)$
is simple (where $Z(G_0)$ denotes the centre of $G_0$).
        \end{result}

We shall recall the details of the groups belonging to the classes
in {\bf (A)}, {\bf (B)} and {\bf (C)} as we need them.

Suppose $\bar K$ is a cap in $PG(r,q)$ such that a subgroup
$\bar{G_0}$ of $P\Gamma L(r+1,q)$ acts transitively on each of $\bar K$
and its complement. Then $\bar{G_0}$ corresponds to a subgroup $G_0$ of
$GL(d,p)$ having three orbits on the vectors of $V(d,p)$, where $p$ is
prime and $p^d=q^{r+1}$. Moreover $G_0$ will contain matrices
corresponding to scalar multiplication by elements of $GF(q)^{*}$. As we
demonstrate shortly, with one exception, $V(d,p)\cdot G_0$ is primitive as
a permutation group, so Liebeck's theorem may be applied. Notice that
since we are interested in groups $G_0$ containing $GF(q)^*$ we avoid the
possibility of two orbits of vectors in $V(d,p)$ giving rise to a single
orbit of points in $PG(r,q)$.

Clearly $G_0$ may be embedded in $\Gamma
L(r+1,q)$. At the beginning of Section 1 of \cite{L}, Liebeck notes that
in his result $G_0\leq GL(d,p)$ is embedded in $\Gamma L(a,p^{d/a})$ with
$a$ minimal. Thus $r+1\geq a$ i.e. $q\leq p^{d/a}$. Moreover in almost all
cases it is clear that the groups he identifies have orbits that are unions of
$1$--dimensional subspaces of $V(a,p^{d/a})$ (excluding the zero vector).
If a $1$-dimensional subspace over $GF(p^{d/a})$ does contains vectors $u,v$ that are linearly independent over $GF(q)$, then $u,v$ and $u+v$ correspond to three collinear points in
$PG(r,q)$ and the orbit in $PG(r,q)$ cannot be a cap. Thus in our setting
we usually have $q=p^{d/a}$: there is just one exception, the class A1, although we have to justify $q=p^{d/a}$ for the class A2.

	\begin{lemma}
Suppose $\bar K$ is a transitive, co-transitive cap in $PG(r,q)$ with $\bar{G_0} \leq P\Gamma L(r+1,q)$ acting transitively on each of $\bar K$ and $PG(r,q) - \bar{K}$ and suppose that $G_0$ is the pre-image of $\bar{G}_0$ in $GL(d,p)$.
Let $H=V(d,p)\cdot G_0$. Then $H$ is imprimitive on $V=V(d,p)$ if and only if $q=2$ and $\bar{K}$ is the complement of a hyperplane.
	\end{lemma}

\Prf Suppose that $H$ is imprimitive on $V$. Let $\Omega$ be a
block containing $0$. Then the two orbits of non--zero vectors of
$G_0$ are $\Omega\setminus{0}$ and $V\setminus\Omega$. Let $u$ and
$v$ be any two vectors in $\Omega$, then $\Omega +v$ is a block
containing $0+v$ and $u+v$ so $\Omega +v=\Omega$. In other words
$u+v$ is in $\Omega$ and so $\Omega$ is a $GF(p)$--subspace of
$V$. More than this $G_0$ contains the scalars in $GF(q)^*$ and so
$\Omega$ is actually a $GF(q)$--subspace. Thus $\Omega$ cannot
correspond to a cap. In $PG(r,q)$ our two orbits consist of points
in a subspace and the complement. A line not in the subspace meets
the subspace in at most one point so the complement cannot form a
cap except perhaps when $p=q=2$ and the subspace has projective
dimension $r-1$. Conversely, as is well known, the complement of a
hyperplane is indeed a cap in $PG(r,2)$ and it is the only way in
which the complement of a subspace is a cap. It is easy to see
that this cap is transitive and co--transitive.
\qed
\\
We recall for the reader that the {\it socle} of a finite group is
the product of its minimal normal subgroups. In our setting
$V(d,p)\cdot G_0$ has $V$ as its unique minimal normal subgroup.

Liebeck's theorem tells us the possibilities for $G_0$ and gives
two orbits of $G_0$ on the non--zero
vectors of $V(d,p)$. We denote these by $K_1$ and $K_2$, and the
corresponding sets of points in $PG(r,q)$ by $\bar{K_1}$ and $\bar
{K_2}$. We assume that neither $K_1$ nor $K_2$ lies in a subspace
of $V(r+1,q)$; given $GF(q)^*\leq G_0$ this means that neither
$K_1$ nor $K_2$ lies in a subspace of $V(d,p)$. We may
henceforth assume that $V(d,p)\cdot G_0$ is a finite primitive affine
permutation group of rank $3$ and degree $p^d$, so we may apply Result \ref{liebeck}.

We begin with the case by case analysis. In many cases we use data from Result \ref{liebeck} and apply Lemmas \ref{bound}, \ref{chordno}, but there are occasions when we need to look at the structure of orbits in detail; there are also occasions when using the structure of the orbits is more illuminating and yet no less efficient than the bound and chord-number arguments.

\section{The infinite classes A}
\subsection{The class A1}

In this case $G_0$ is a subgroup of $\Gamma L(1,p^d)$ containing $GF(q)^*$. Such a subgroup 
is generated by $\omega^N$ and $\omega^e\alpha^s$, for some $N, e, s$ where $\omega$ is a primitive element of $GF(p^d)$ and $\alpha$ is the generating automorphism $x \mapsto x^p$ of $GF(p^d)$; if we write $p^d=q^a$, then $N$ divides $(q^a-1)/(q-1)$. Let $H_0$ be the subgroup of $\Gamma L(1,p^d)$ generated by $\omega^N$. Then $H_0$ is a Singer subgroup of $GL(1,p^d)$ and the orbits of $H_0$ in $PG(r,q)$ are called Singer orbits. Clearly if $G_0$ has two orbits in $PG(r,q)$, then each orbit is the union of Singer orbits. If the smaller orbit is to be a cap, then each Singer orbit must itself be a cap. A precise criterion for deciding when Singer orbits are caps in $PG(r,q)$ is given by Sz\H onyi \cite[Proposition 1]{TS}.

Precise criteria for there to be two orbits for $G_0$ on non--zero vectors of $V(d,p)$ are given by Foulser and Kallaher \cite{FK}. These involve numbers $m_1$ and $v$ such that the primes of $m_1$ divide $p^s-1$, $v$ is a prime $\neq 2$ and ord$_vp^{sm_1}=v-1$, $(e,m_1)=1$, $m_1s(v-1)\vert d$, $N=vm_1$. The orbit lengths are $m_1(p^d-1)/N$ and $(v-1)m_1(p^d-1)/N$. Notice that when $p=2$ the smaller orbit has odd size. Hill \cite{RH2} suggests the possibility of transitive, co--transitive caps of size $78$ in $PG(5,4)$ and $430$ in $PG(6,4)$. It is now clear that these cannot be caps from class $A1$ and our main theorem then shows that they cannot be caps at all.

\subsection{The class A2}

$G_0$ preserves a direct sum $V_1 \oplus V_2$, where $V_1,V_2$ are subspaces of $V(d,p)$. One orbit must be $K_1=(V_1 \cup V_2)-\{0 \}$ and the other $K_2= \{v_1 +v_2: 0\neq v_1 \in
V_1, 0\neq v_2 \in V_2 \}$. We first show that $V_1,V_2$ are subspaces over $GF(q)$. Observe that for any $\lambda \in GF(q)^* \leq G_0$, $\lambda V_1=V_1$ or $V_2$ and let $F=\{\lambda \in GF(q)^*: \lambda V_1=V_1\} \cap \{0\}$. Then $F$ is a subfield of $GF(q)$ having order greater than $q/2$ so must be $GF(q)$. It is now clear that $V_1,V_2$ are subspaces of $V(r+1,q)$ of dimension $t=(r+1)/2$. Given that $r \geq 2$, we must have $t \geq2$, so $\bar{K_1}$ contains lines of $PG(r,q)$ and cannot be a cap. Moreover $| \bar{K_1} |=2(q^t-1)/(q-1) <(q^{r+1}-1)/2$ so $\bar{K_1}$ is the smaller orbit and therefore $\bar{K_2}$ cannot be a cap.

\subsection{The class A3}

$G_0$ preserves a tensor product $V_1 \otimes V_2$ over $GF(q)$, with $V_1$ having dimension 2 over $GF(q)$. One orbit must be $K_1=\{v_1 \otimes v_2: 0\neq v_1 \in V_1, 0\neq v_2 \in V_2
\}$ and the other $K_2=V-(K_1 \cup \{0\}) $.

Consider the $GF(q)$--subspace $V_1 \otimes v_2$
of V for some $0 \neq
v_2 \in V_2$. It has dimension $2$ in $V(r+1,q)$ so corresponds to
a line in $PG(r,q)$ inside $\bar{K_1}$. Hence $\bar{K_1}$ is not a cap.

Let $b$ be the dimension of $V_2$ over $GF(q)$. Then $r+1=2b$ and $|\bar{K_1}|=(q+1)(q^b-1)/(q-1)$ (\cite[Table12]{L}) so $|\bar{K_2}|=q(q^b-1)(q^{b-1}-1)/(q-1)>|\bar{K_1}|$ except when $q=2,b=2$ (i.e., $r+1=d=4$). Thus there is only one case in which $\bar{K_2}$ can possibly be a cap.

Suppose that $q=p=2$ and $d=4$, i.e. we are reduced to considering caps in $PG(3,2)$. In $PG(3,2)$, we see that $\vert\bar{K_1}\vert = 9$ and $\vert\bar{K_2}\vert = 6$ . Thus here $\bar{K_1}$ is too big and for $\bar{K_2}$ it is simplest to note that $(6.5.1)/(2.9) \notin \mathbb{Z}$, so neither is a cap (by Lemmas \ref{bound} and  \ref{chordno}).

\subsection{The class A4}

$G_0 \unrhd SL(a,s)$ and $p^d = s^{2a}$. Here $q = s^2$ and $p^d =
q^a$ with $SL(a,s)$ embedded in $GL(d,p)$ as a subgroup of
$SL(a,q)$: let $e_1, e_2, ...,e_a$ be a basis for $V$ over $GF(q)$
then with respect to this basis $SL(a,s)$ consists of the matrices in
$SL(a,q)$ having every entry in $GF(s)$. If $G_0$ has two orbits on
non-zero vectors of $V$ then those orbits must be $K_1 =\{\gamma \sum
\lambda_i e_i$ ($\lambda_i \in GF(s)$, not all $0$),$0 \neq \gamma \in
GF(q)\}$ and $K_2$ the set of all remaining non-zero vectors. In other
words $\bar{G_0}$ preserves a subgeometry of $PG(r,q)$. We have $r>1$ so
that $a\geq 3$. Thus three collinear points of $PG(r,s)$ are still three
collinear points in $PG(r,q)$ and so $\bar{K_1}$ is not a cap.

Let us turn to $\bar{K_2}$. As noted above, $r>1$ so $a \geq 3$. Let
$u=e_1 + \sigma e_2, v=e_2 + \sigma e_3$, where $\sigma\in GF(q)\setminus
GF(s)$. Then $u,v$ and $u+v =e_1 + (\sigma +1)e_2 + \sigma e_3 \in K_2$
correspond to collinear points of $PG(r,q)$, all in $\bar{K_2}$. Hence
$\bar{K_2}$ is not a cap.

\subsection{The class A5}

$G_0 \unrhd SL(2,s)$ and $p^d = s^6$. Here $q = s^3$ and $p^d =
q^2$ with $SL(2,s)$ embedded in $GL(d,p)$ as a subgroup of
$SL(2,q)$. However $r=1$ in this case so it does not concern us.

\subsection{The class A6}

$G_0 \unrhd SU(a,q')$ and $p^d = ((q')^2)^a$. In this case $q={(q')}^2$.
Here one orbit $K_1$ consists of the non-zero isotropic vectors and the
other orbit $K_2$ consists of the non-isotropic vectors with respect to an
appropriate non-degenerate hermitian form. Each orbit is a union of
$1$--dimensional subspaces of $V(a,q)$ (excluding the zero vector). To
begin with, a non--isotropic line of $PG(r,q)$ contains at least three
isotropic points, i.e., three points of $\bar{K_1}$. Therefore $\bar{K_1}$
cannot be a cap.

Now consider $\bar{K_2}$. Given $a \geq 3$, consider a line of
$PG(r,q)$ that is isotropic but not totally isotropic, then it
contains one point of $\bar{K_1}$ and $q\geq 4$ points of
$\bar{K_2}$. Hence $\bar{K_2}$ is not a cap.

\subsection{The class A7}
$G_0 \unrhd \Omega^{\pm}(a,q)$ and $p^d = (q)^{a}$ with $a$ even (and if
$q$ is odd , $G_0$ contains an automorphism interchanging the two orbits
of $\Omega^{\pm}(a,q)$ on non-singular 1-spaces). The arguments here are
similar to the Unitary case. $K_1$ consists of the non-zero singular
vectors and $K_2$ consists of the non-singular vectors. Let $b$ be the
Witt index of the approppriate quadratic form on $V(a,q)$ i.e., the
dimension of a maximal totally singular subspace.
 Then $a$ is one of $2b,2b+2$.
Any totally singular line would be a line of $PG(r,q)$ lying
inside $\bar{K_1}$. Given that $a\geq 3$, it follows
that the only possibility for $\bar{K_1}$ being a cap is when
$\bar{K_1}$ is an elliptic quadric in $PG(3,q)$. In passing we note that for odd $q$, the necessary automorphism is contained in $G_0$ only when $q$ is square; in this case and in the case $q$ even, the elliptic quadric gives a well known cap.

Let us turn to $\bar{K_2}$. Any anisotropic line of $PG(r,q)$ lies
inside $\bar{K_2}$ so $\bar{K_2}$ can never be a cap.

\subsection{The class A8}

$G_0 \unrhd SL(5,q)$ and $p^d = (q)^{10}$ (from the action of
$SL(5,q)$ on the skew square $\bigwedge ^2 (V(5,q))$. Here one orbit
of non-zero vectors must be $K_1=\{0 \neq u \bigwedge v: u,v \in
V(5,q)\}$ with the other non-zero vectors belonging to $K_2$. One can argue in a similar manner to the case of the tensor product. However it is quicker here to note that the orbits of $\bar{G_0}$ on $PG(r,q)$ have sizes $k=(q^5-1)(q^2+1)/(q-1)$ and $m=q^2(q^5-1)(q^3-1)/(q-1)$ (\cite[Table12]{L}) with $k < m$ for all values of $q$. The chord-number is then given by $c=k(k-1)(q-1)/2m$ by Lemma \ref{chordno} i.e., $c=(q^2+1)(q^3+q+1)/2q \notin \mathbb{Z}$. Hence neither $\bar{K_1}$ nor $\bar{K_2}$ is a cap.

\subsection{The class A9}
$G_0/Z(G_0)\unrhd\Omega(7,q)\cdot Z_{(2,q-1)}$ and $p^d=q^8$ (from the
action of $B_3(q)$ on a spin module) \cite{CC}, \cite{KL}. The study of
Clifford algebras leads to the construction of "spin modules" for
$P\Omega (m,q)$. When $m=8$ this leads to the triality automorphism of
$P\Omega^+ (8,q)$. One finds that it is possible (via this automorphism)
to embed $\Omega (7,q) \cong P\Omega (7,q)$ inside $P\Omega^+ (8,q)$ as
an irrdeucible subgroup. The important thing from our point of view is
that two non--trivial orbits of $G_0$ must be the set of all non--zero
singular vectors and the set of all non--singular vectors with respect
to a non--degenerate quadratic form on $V(8,q)$. In this setting the
arguments employed for class $A7$ apply: neither orbit can be a cap.

\subsection{The class A10}

$G_0/Z(G_0)\unrhd P\Omega^+(10,q)$ and $p^d=q^{16}$ (from the action of $D_5(q)$ on a spin module) \cite{CC}, \cite{KL}. Once again we have a spin representation, this time of
$P\Omega^+(10,q)$ on $PG(15,q)$. On this occasion it is quickest to work from the orbit lengths.

The orbits of $\bar{G_0}$ on $PG(r,q)$ have sizes $k=(q^8-1)(q^3+1)/(q-1)$ and $m=q^3(q^8-1)(q^5-1)/(q-1)$ (\cite[Table12]{L}) with $k < m$ for all values of $q$. The chord-number is then given by $c=k(k-1)(q-1)/2m$ by Lemma \ref{chordno} i.e., $c=(q^3+1)(q^5+q^2+1)/2q^2 \notin \mathbb{Z}$. Hence neither $\bar{K_1}$ nor $\bar{K_2}$ is a cap.

\subsection{The class A11}

$G_0 \unrhd Sz(q)$ and $p^d = (q)^4$, with $q\geq 8$ an odd power of $2$
(from the embedding $Sz(q) \leq Sp(4,q)$). Here the smaller orbit
$\bar{K_1}$ on $PG(3,q)$ is a Suzuki--Tits ovoid containing $q^2+1$ points
and this is indeed a cap \cite{Ti}, \cite[16.4.5]{H}.

\section{The Extraspecial classes}

In most cases here $G_0\leq M$ where $M$ is the normalizer in
$\Gamma L(a,q)$ of a $2$--group $R$, where $p^d=(q)^a$ and
$a=2^m$ for some $m\geq 1$; either $R$ is an extraspecial group
$2^{1+2m}$ or $R$ is isomorphic to $Z_4\circ 2^{1+2m}$. In all
cases here $p$ is odd. There are two types of extraspecial group
$2^{1+2m}$, denoted $R_1^m$ and $R_2^m$; the first of these
has the structure $D_8\circ D_8\circ\dots D_8$ ($m$ copies)
and the second $D_8\circ D_8\circ\dots\circ D_8\circ Q_8$ ($m-1$ copies of
$D_8$), where $D_8$ and $Q_8$ are respectively the dihedral and quaternion
groups of order $8$, and '$\circ$' indicates a central product. The group
$Z_4\circ 2^{1+2m}$ is again a central product, this time $Z_4\circ
D_8\circ D_8\circ\dots\circ D_8$ ($m$ copies of $D_8$) and is denoted by
$R^m_3$. Notice that $R$ modulo its centre is an elementary abelian
$2$--group, i.e. a $2m$--dimensional vector space over $GF(2)$ and in fact
$M/RZ$ ($Z$ being the centre of $\Gamma L(a,q)$) may be embedded in
$GSp(2m,2)$. In just one case $G_0\leq M$ with $M$ the normalizer in
$\Gamma L(3,4)$ of a $3$--group of order $27$. We record from \cite[Table
13]{L} that in this case the non--trivial orbit sizes of $G_0$ on $V(3,4)$
are $27$ and $36$, i.e. the point orbit sizes in $PG(2,4)$ are $9$ and
$12$, but the largest possible size of a cap (here better termed an arc)
in $PG(2,4)$ is $6$. Hence there are no caps here and we
may henceforth assume that $R$ is a $2$--group, with $p$ odd.

There are sixteen instances where $G_0$ has two non--trivial
orbits on $V(d,p)\simeq V(a,q)$, but ten of these have $a=2$
(i.e. $m=1$) and so refer to action on a projective line, i.e.
$r<2$; note that two of these cases have $q>p$.
Thus we concentrate on the remaining six cases.
In each of these cases $q=p$ and in all but the last case the
vector space is $V(4,p)$. In the last case the vector space is $V(8,3)$. Four cases follw immediately from known bounds - they are listed in the table below.
\vspace{.5cm}

\begin{center}
\begin{tabular}{| c | c | c | c | c |}  \hline
p=q & r & R & smaller orbit size & max. cap size\\ \hline
3&3& $R_1^2$ & 16 & 10 \\  \hline
5&3&$R_2^2$ & 60 & 26 \\ \hline
5&3&$R_3^2$ & 60 & 26 \\ \hline
7&3&$R_2^2$ & 80 & 50 \\ \hline

\end{tabular}
\end{center}

\vspace{.5cm}

\textbf{The case $\mathbf{p=q=3}$, $\mathbf{r=7}$, $\mathbf{R=R_2^3}$.}

In this case smaller orbit of $\bar{G_0}$ on $PG(7,3)$ has size $720$, while the maximum
size for a cap in $PG(7,3)$ is only known to be $\leq 729$. Instead we use Lemma \ref{chordno}: the larger orbit has size $2560$ and $(720.719.2)/(2.2560) \notin \mathbb{Z}$.\\

\textbf{The case $\mathbf{p=q=3}$, $\mathbf{r=3}$, $\mathbf{R=R_2^2}$.}

Here Liebeck notes that $R$ has five orbits of size $16$ on
$V(4,3)$ and $M$ permutes these orbits acting as $S_5$, the
symmetric group of degree $5$. Thus there are a number of
possibilities for $G_0$ having two non--trivial orbits on
$V(4,3)$. However it is straightforward to construct generating
matrices for $R$ and we see immediately that one orbit of size
$16$ on $V(4,3)$ cannot correspond to a cap in $PG(3,3)$.
Therefore none of the orbits of size $16$ can correspond
to a cap and hence no possible choices of $G_0$ can give rise to a
cap.

\section{The Exceptional classes}

Finally we turn to the exceptional classes where the socle $L$ of
$G_0/Z(G_0)$ is simple. There are just thirteen different
possibilities for $L$, although on occasion more than one
possibility for $G_0$ corresponds to a given $L$. For example for
$L=A_5$ there are seven different possibilities for $G_0$ (one of which leads to a single orbit in $PG(d-1,p)$); however all of these lead to $r<2$ and so do not concern us.

We employ a variety of techniques to tackle these cases. Liebeck
\cite[Table 14]{L} gives the orbit sizes in $V(d,p)$ and sometimes
we can use these to rule out the possibility of caps. On other
occasions we can use the fact that the chord-number is an integer. On two occasions, neither of these appraoches works and we have to investigate the known structure of the smaller orbit. There remain two cases where a cap does occur.\\

\textbf{The cases where caps occur.}

When $L=A_6$ and $(d,p)=(6,2)$, $L$ admits an embedding in
$PSL(3,4)$ (so here $q=4$) and $G_0$ has an orbit of size $6$. In
fact this in a hyperoval in $PG(2,4)$ \cite{B},\cite{G} so we do have a
cap.

When $L=M_{11}$ and $(d,p)=(5,3)$ there is a representation of $L$ in which one orbit has size $11$ and in fact this is a cap. In passing we note that this cap arises as an orbit of a Singer cyclic subgroup of $PG(4,3)$ \cite{CK}; moreover $PG(4,3)$ is partitioned into eleven
$11$--caps (the eleven orbits of the Singer cyclic subgroup). Note also that there is a second representation of $L=M_{11}$ on $PG(4,3)$ (see below). In fact both representations appear in the context of the ternary Golay code \cite[Ch. 6]{Asch2}.\\

\textbf{Cases where known bounds rule out caps.}

In each of the following cases the smaller orbit is larger than the known upper bound for a cap size, so cannot be a cap. In the table $k$ is the smaller orbit size.

\vspace{.5cm}

\begin{center}
\begin{tabular}{| c | c | c | c|c | c |}  \hline
$L$ & $(d,p)$ & $r$ & $q$ & $k$ & max. cap size\\ \hline 
$A_6$&$(4,5)$& 3&5 & 36 & 26\\  \hline
$A_7$&$(4,7)$& 3&7 & 120 & 50\\  \hline
$M_{11}$&$(5,3)$& 4&3 & 55 & $\leq 27$\\  \hline
$J_2$&$(6,5)$& 5&5 & 1890 & $\leq 625$\\  \hline
$J_2$&$(12,2)$& 5&4 & 525 & $\leq 256$\\  \hline

\end{tabular}
\end{center}

\vspace{.5cm}

\textbf{Cases where $c$ an integer rules out caps.}

In each of the following cases a calculation $c={k(k-1)(q-1)/2m}$ yields a non-integer and so by Lemma \ref{chordno}, the smaller orbit does not correspond to a cap. In the table $k$ is the smaller orbit size and $m$ the larger orbit size.

\vspace{.5cm}

\begin{center}
\begin{tabular}{| c | c | c | c|c | c |}  \hline
$L$ & $(d,p)$ & $r$ & $q$ & $k$ & $m$\\ \hline 
$A_9$&$(8,2)$& 7&2 & 120 & 135\\  \hline
$A_{10}$&$(8,2)$&7 &2 & 45 & 210\\  \hline
$L_2(17)$&$(8,2)$& 7&2 & 102 & 153\\  \hline
$M_{24}$&$(11,2)$& 10&2 & 276 & 1771\\  \hline
$M_{24}$&$(11,2)$& 10&2 & 759 & 1288\\  \hline
$Suz$ or $J_4$&$(12,3)$& 11& 2 & 65520 & 465920\\  \hline

\end{tabular}
\end{center}

\vspace{.5cm}

\textbf{The case $\mathbf{L=A_7}$ and $\mathbf{(d,p)=(8,2)}$.}

Here $L$ is embedded in $PSL(4,4)$ (so $q=4$). In fact $L$ may
actually embedded in $A_8\simeq PSL(4,2)\leq PSL(4,4)$. The group
$A_8$ and therefore $A_7$ preserve a subgeometry whose $15$ points
form the smaller orbit. There are numerous examples of three
points on a line in the subgeometry. Thus we have no caps.

\textbf{The case $\mathbf{L=PSU(4,2)}$ and $\mathbf{(d,p)=(4,7)}$.}

The vectors in the smaller orbit are given by Liebeck \cite[Lemma
3.4]{L}:
$$
(\theta;0,0,0),\quad (0;\theta,0,0),\quad
(0;\omega^a,\omega^b,\omega^c),\quad
(\omega^a;0,\omega^b,-\omega^c),
$$
(together with all scalar multiples) where $\theta =\omega =2$;
$a,b,c$ take any integral values; and the last three coordinates
may be permuted cyclically. It suffices here to observe that
$(1;0,0,0)$, $(1;0,1,6)$ and $(2;0,1,6)$ all lie in this orbit and
give three collinear points in $PG(3,7)$. So no cap arises here.

\bigskip

\noindent {\it AMS Mathematics Subject Classification: Primary 51E22; Secondary 20B15, 20B25

\noindent Keywords and Phrases: Caps, Rank 3 permutation groups}
\vspace{0.5cm}

\noindent
{\bf A. Cossidente},
Dipartimento di Matematica,
Universit\`a della Basilicata,
via N.Sauro 85, 85100 Potenza, Italy.\\
e-mail: cossidente@unibas.it
\par\noindent
{\bf O.H. King},
Department of Mathematics,
The University of Newcastle,
Newcastle Upon Tyne,
NE1 7RU, United Kingdom.\\
e-mail: o.h.kink@ncl.ac.uk
\end{document}